\input amstex\documentstyle{amsppt}  
\pagewidth{12.5cm}\pageheight{19cm}\magnification\magstep1  
\topmatter
\title Positive structures in Lie theory\endtitle
\author G. Lusztig\endauthor
\address{Department of Mathematics, M.I.T., Cambridge, MA 02139}\endaddress
\thanks{Supported by NSF grant DMS-1566618.}\endthanks
\endtopmatter   
\document

\define\uI{\un I}

\define\ui{\un i}
\define\uj{\un j}

\define\bQ{\bar Q}

\define\si{\sim}

\define\sqc{\sqcup}

\define\part{\partial}
\define\emp{\emptyset}

\define\m{\mapsto}
\define\do{\dots}

\define\lra{\leftrightarrow}

\define\T{\times}

\define\nl{\newline}
\redefine\i{^{-1}}

\define\un{\underline}

\define\a{\alpha}

\redefine\c{\chi}

\define\e{\epsilon}
\define\et{\eta}

\redefine\o{\omega}

\define\ph{\phi}
\define\ps{\psi}

\redefine\t{\tau}

\define\hh{\bold h}
\define\ii{\bold i}

\define\kk{\bold k}

\define\BB{\bold B}
\define\CC{\bold C}

\define\NN{\bold N}

\define\RR{\bold R}

\define\ZZ{\bold Z}

\define\cb{\Cal B}

\define\cj{\Cal J}

\define\co{\Cal O}
\define\cp{\Cal P}

\define\bP{\bar P}

This paper is based on a lecture given at the International Consortium of Chinese Mathematicians,
Taipei, December 2018.

\subhead 0.1\endsubhead
In late 19th century and early 20th century, a new branch of mathematics was born: Lie theory or the 
study of Lie groups and Lie algebras (Lie, Killing, E.Cartan, H.Weyl). It has become a central part 
of mathematics with applications everywhere. More recent developments in Lie theory are as follows.

-Analogues of simple Lie groups over any field (including finite
fields where they explain most of the finite simple groups): Chevalley 1955;

-infinite dimensional versions of the simple Lie algebras/simple Lie groups: Kac and Moody 1967, 
Moody and Teo 1972;

-theory of quantum groups: Drinfeld and Jimbo 1985.

\subhead 0.2\endsubhead
In Lie theory to any Cartan matrix one can associate a simply connected Lie group $G(\CC)$;
Chevalley replaces $\CC$ by any field $\kk$ and gets a group $G(\kk)$. In \cite{L94} we have defined 
the totally positive (TP) submonoid $G(\RR_{>0})$ of $G(\RR)$ and its ``upper triangular'' part
$U^+(\RR_{>0})$.  In this lecture we will review the TP-monoids $G(\RR_{>0})$, $U^+(\RR_{>0})$
attached to a 
Cartan matrix, which for simplicity is assumed to be 
simply-laced. 
In \cite{L94} the nonsimply laced case is treated by reduction to the simply laced case. 

\subhead 0.3\endsubhead
The total positivity theory in \cite{L94} was a starting point for 

-a solution of Arnold's problem for real flag manifolds, Rietsch 1997;

-the theory of cluster algebras, Fomin, Zelevinsky 2002;

-a theory of TP for the wonderful compactifications, He 2004;

-higher Teichm\"uller theory, Fock, Goncharov 2006;

-the use of the TP grassmannian in physics, Postnikov 2007, Arkani-Hamed, Trnka 2014;

-a theory of TP for the loop group of $GL_n$, Lam, Pylyavskyy 2012;

-a theory of TP for certain nonsplit real Lie groups, Guichard-Wienhard 2018;

-a new approach to certain aspects of quantum groups, Goncharov, Shen.

\subhead 0.4\endsubhead
Schoenberg (1930) and Gantmacher-Krein (1935) (after initial contributions of Fekete and Polya (1912))
defined the notion of TP matrix in $GL_n(\RR)$: a matrix all of whose $s\T s$ minors are $\ge0$ for any $s$.
Gantmacher and Krein showed that if for any $s$, all $s\T s$ minors of a matrix $A$ are $>0$ then the 
eigenvalues of $A$ are real, distinct and $>0$. 
For example, the Vandermonde matrix $(A_{ij})$, $A_{ij}=x_i^{j-1}$ with $x_1<x_2<\do<x_n$ real and $>0$ is of 
this type. According to Polya and Szeg\"o, the matrix $(A_{ij})$, $A_{ij}=\exp(x_iy_j)$ with 
$x_1>x_2>\do>x_n$, $y_1>y_2>\do>y_n$ real is also of this type.

The TP matrices in $GL_n(\RR)$ form a monoid under multiplication. This monoid is 
generated by diagonal matrices with $>0$ entries on diagonal and by matrices which have $1$ on 
diagonal and a single nonzero entry off diagonal which is $>0$ (Whitney, Loewner, 1950's). 
Our definition \cite{L94} of the TP part of any $G(\RR)$ was inspired by the work of Whitney, Loewner.

However, to prove properties of the resulting
monoid (such as the generalization of the Gantmacher-Krein theorem) I had 
to use the canonical bases in quantum groups (discovered in \cite{L90}) and their positivity properties. 
The role of $s\T s$ minors is played in the general case by the canonical basis of \cite{L90}. 
Unlike in \cite{L94}, here we define $G(\RR_{>0})$ by generators and relations, independently of 
$G(\RR)$. Surprisingly, this definition of $G(\RR_{>0})$ is simpler than that of $G(\RR)$ (see 
\cite{ST}). From it 
one can recover the Chevalley groups $G(\kk)$ for any field $\kk$. Namely, the relations
between the generators of $G(\RR_{>0})$ involve only rational functions with integer coefficients. 
They make sense over $\kk$ and they give rise to a ``birational form'' of a
semisimple group over $\kk$. 
This is the quotient field of the coordinate ring of $G(\kk)$; then 
$G(\kk)$ itself appears as a subgroup of the automorphism group of 
this field. In this approach the form $G(\RR_{>0})$ is the most basic, 
all other forms are deduced from it. 

\subhead 0.5\endsubhead
We now describe the content of various sections.
In \S1 we define a positive structure on a set.  Such structures have appeared in \cite{L90}, \cite{L94}
in connection with various objects in Lie theory. In \S2 we define the monoid $U^+(\RR_{>0})$.
In \S3 we define the monoid $G(\RR_{>0})$.
In \S4 we use this monoid to recover the Chevalley groups over a field.
In \S5 we define the non-negative part of a flag manifold.

\head 1. Positive structures\endhead
\subhead 1.1\endsubhead
The TP monoid can be defined not only over $\RR_{>0}$ but over a structure $K$ in which 
addition, multiplication, division (but no substraction) are defined. In \cite{L94} three such 
$K$ were considered.

(i) $K=\RR_{>0}$; 

(ii) $K=\RR(t)_{>0}$, the set of $f\in\RR(t)$ of form $f=t^ef_0/f_1$ for some 

$f_0,f_1$ in $\RR[t]$ with constant term in $\RR_{>0},e\in\ZZ$ ($t$ is an 

indeterminate); 

(iii) $K=\ZZ$.
\nl
In case (i) and (ii), $K$ is contained in a field $\RR$ or $\RR(t)$ and the sum and product 
are induced from that field. In case (iii) we consider a 
new sum $(a,b)\m\min(a,b)$ and a new product $(a,b)\m a+b$. A 4th case is 

(iv) $K=\{1\}$
\nl
with $1+1=1,1\T1=1$.

In each case $K$ is a semifield (a terminology of Berenstein, Fomin, Zelevinsky 1996): a 
set with two operations, $+$, $\T$, which is an abelian group with respect to $\T$, an abelian 
semigroup with respect to $+$ and in which $(a+b)c=ac+bc$ for all $a,b,c$.  We fix a semifield
$K$. 
There is an obvious 
semifield homomorphism $K\to\{1\}$. We denote by $(1)$ the unit element of $K$ with respect to $\T$.

\subhead 1.2\endsubhead
In \cite{L94} we observed that there is a semifield homomorphism $\a:\RR(t)_{>0}\to\ZZ$ given by 
$t^ef_0/f_1\m e$ which connects geometrical objects over $\RR(t)_{>0}$) with piecewise linear objects 
involving only integers. I believe that this was the first time that such a connection (today known 
as tropicalization) was used in relation to Lie theory. 

\subhead 1.3\endsubhead
For any $m\in\ZZ_{>0}$ let $\cp_m$ be set of all nonzero polynomials in $m$ indeterminates $X_1,X_2,\do,X_m$
with coefficients in $\NN$. 

A function $(a_1,a_2,\do,a_m)\m(a'_1,a'_2,\do,a'_m)$ from $K^m$ to $K^m$ is said to be {\it admissible}
if for any $s$ we have $a'_s=P_s(a_1,a_2,\do,a_m)/Q_s(a_1,a_2,\do,a_m)$ where $P_s,Q_s$ are in $\cp_m$.
(This ratio makes sense since $K$ is a semifield.) In the case where $K=\ZZ$, such a function is 
piecewise-linear. If $m=0$, the unique map $K^0@>>>K^0$ is considered to be admissible ($K^0$ is a point.)

\subhead 1.4\endsubhead
A {\it positive structure} on a set $X$ consists of a family of bijections $f_j:K^m@>\si>>X$ 
(with $m\ge0$ fixed) indexed by $j$ in a finite set $\cj$, such that $f_{j'}\i f_j:K^m@>>>K^m$ is 
admissible for any $j,j'$ in $\cj$; the bijections $f_j$ are said to be the {\it coordinate charts} of 
the positive structure. The results of \cite{L94}, \cite{L97}, \cite{L98}, can be interpreted 
as saying that various objects in Lie theory admit natural positive structures.

\head 2. The monoid $U^+(K)$\endhead
\subhead 2.1. The Cartan matrix\endsubhead
We fix a finite graph; it is a pair consisting of two finite sets $I,H$ and a map which to each
$h\in H$ associates a two-element subset $[h]$ of $I$.
The Cartan matrix  $A=(i:j)_{i,j\in I}$ is given by $i:i=2$ for all $i\in I$ while if
$i,j$ in $I$ are distinct then $i:j$ is $-1$ times the number of $h\in H$ such that
$[h]=\{i,j\}$. 

An example of a Cartan matrix is:

$I=\{i,j\}$, $A=\left(\matrix 2&-1\\-1&2\endmatrix\right)$.
\nl
We fix a Cartan matrix $A$. For applications to Lie theory $A$ is assumed to be 
positive definite. But several of the subsequent definitions make sense without this assumption.

We attach to $A$ and a field $\kk$ a group $G(\kk)$. When $A$ is positive definite, $G(\kk)$
is the group of $\kk$-points of a simply connected semisimple split algebraic group of type $A$ over $\kk$.
Without the assumption that $A$ is positive definite, the analogous group $G(\kk)$ (with $\kk$ of
characteristic $0$) has been defined
 in \cite{MT}, \cite{Ma},\cite{Ti}.

We will associate to $A$ and $K$ a monoid $G(K)$ and a submonoid $U^+(K)$ of $G(K)$. 
In the case where $K=\RR_{>0}$ (resp. $K=\RR(t)_{>0}$), $G(K)$ and $U^+(K)$ can be viewed as
submonoids of $G(\kk)$ where $\kk=\RR$ (resp. $\kk=\RR(t)$). In the case where $K=\RR_{>0},\kk=\RR,
G(\RR)=SL_n(\RR)$, $U^+(K)$ is the monoid of TP matrices in $G(\RR)$ which are 
upper triangular with $1$ on diagonal. We first define $U^+(K)$. 

\subhead 2.2\endsubhead
Let $U^+(K)$ be the monoid (with $1$) with generators 
$i^a$ with $i\in I$, $a\in K$ and relations 

$i^ai^b=i^{a+b}$ for $i\in I$, $a,b$ in $K$;

$i^aj^bi^c=j^{bc/(a+c)}i^{a+c}j^{ab/(a+c)}$ for $i,j\in I$ with $i:j=-1$, $a,b,c$ in $K$;

$i^aj^b=j^bi^a$ for $i,j\in I$ with $i:j=0$, $a,b$ in $K$.
\nl
(In the case where $K=\ZZ$, relations of the type considered above 
involve piecewise-linear functions; they first appeared in \cite{L90} in 
connection with the parametrization of the canonical basis.)
The definition of $U^+(K)$ is reminiscent of the definition of the Coxeter group attached to $A$.
In the case where $K=\ZZ$ and $A$ is positive definite the definition of $U^+(K)$ given above
first appeared in \cite{L94, 9.11}.

\subhead 2.3\endsubhead
When $A=\left(\matrix 2&-1\\-1&2\endmatrix\right)$, $K=\RR_{>0}$, we can 
identify $U^+(K)$ with the submonoid of $SL_3(\RR)$ generated by

$\left(\matrix 1&a&0\\0&1&0\\0&0&1\endmatrix\right)$, 
$\left(\matrix 1&0&0\\0&1&b\\0&0&1\endmatrix\right)$, 

with $a,b$ in $\RR_{>0}$.

\subhead 2.4\endsubhead
Let $W$ be the Coxeter group attached to $A$. It has generators $i$ with $i\in I$ and relations
$ii=1$ for $i\in I$; $iji=jij$ for $i,j\in I$, $i:j=-1$; $ij=ji$ for $i,j\in I$, $i:j=0$.
Let $\co_w$ be the set of reduced expressions of $w$ that is the set of 
sequences $(i_1,i_2,\do,i_m)$ in $I$ such that $i_1i_2\do i_m=w$ in $U^+(\{1\})$
where $m$ is minimum. We write $m=|w|$ ($=$length of $w$). 

When $K=\{1\}$, $U^+(K)$ is the monoid (with $1$) with generators $i^1$ 
with $i\in I$ and relations $i^1i^1=i^1$ for $i\in I$; $i^1j^1i^1=j^1i^1j^1$ for 
$i,j\in I$, $i:j=-1$; $i^1j^1=j^1i^1$ for $i,j\in I$, $i:j=0$.
By a lemma of Iwahori and Matsumoto we have can identify (as sets) $W=U^+(\{1\})$ by
$w=i_1i_2\do i_m\lra i_1^1i_2^1\do i_m^1$ for any 
$(i_1,i_2,\do,i_m)\in\co_w$. This bijection is not compatible with the monoid structures.

\subhead 2.5\endsubhead
The semifield homomorphism $K\to\{1\}$ induces a map of monoids
$U^+(K)\to U^+(\{1\})$. Let $U^+_w(K)$ be the fibre over $w\in U^+(\{1\})$. 
We have $U^+(K)=\sqc_{w\in W}U^+_w(K)$. 

We now fix $w\in W$. It turns out that the set $U^+_w(K)$ can be parametrized by $K^m$, in 
fact in many ways, indexed by $\co_w$. For $\ii=(i_1,i_2,\do,i_m)\in\co_w$ 
we define $\ph_\ii:K^m\to U^+_w(K)$ by 

$\ph_\ii(a_1,a_2,\do,a_m)=i_1^{a_1}i_2^{a_2}\do i_m^{a_m}$. 
\nl
This is a bijection. Now $U^+_w(K)$ together with the bijections $\ph_\ii:K^m\to U_w^+(K)$ is an example 
of a positive structure. (We will see later other such structures.)

\subhead 2.6\endsubhead
Let $w\in W,m=|w|$. In the case $K=\ZZ$, $U^+_w(\NN):=\ph_\ii(\NN^m)\subset U^+_w(\ZZ)$
is independent of $\ii\in\co_w$. We set $U^+(\NN)=\sqc_{w\in W}U^+_w(\NN)$; this is a subset of $U^+(\ZZ)$.

When $W$ is finite, let $w_I$ be
the element of maximal length of $W$. Let $\nu=|w_I|$. Now $U^+_{w_I}(\NN)$ 
was interpreted in \cite{L90} as an indexing set for the canonical basis 
of the plus part of a quantized enveloping algebra. A similar 
interpretation holds for $U^+_w(\NN)$ when $W$ is not necessarily finite 
and $w$ is arbitrary, using \cite{L96, 8.2}.

\head 3. The monoid $G(K)$\endhead
\subhead 3.1\endsubhead
In order to define the monoid $G(K)$ we consider besides $I$, two 
other copies $-I=\{-i;i\in I\}$, $\uI=\{\ui;i\in I\}$ of $I$, in obvious 
bijection with $I$. For $\e=\pm1$, $i\in I$ we write $\e i=i$ if $\e=1$, 
$\e i=-i$ if $\e=-1$. 

Let $G(K)$ be the monoid (with $1$) with generators $i^a,(-i)^a,\ui^a$ 
with $i\in I,a\in K$ and the relations below.

(i) $(\e i)^a(\e i)^b=(\e i)^{a+b}$ for $i\in I$, $\e=\pm1$, $a,b$ in $K$;

(ii) $(\e i)^a(\e j)^b(\e i)^c=(\e j)^{bc/(a+c)}(\e i)^{a+c}(\e j)^{ab/(a+c)}$ 

for $i,j$ in $I$ with $i:j=-1$, $\e=\pm1$, $a,b,c$ in $K$;

(iii) $(\e i)^a(\e j)^b=(\e j)^b(\e i)^a$ 

for $i,j$ in $I$ with $i:j=0$, $\e=\pm1$, $a,b$ in $K$;

(iv) $(\e i)^a(-\e i)^b=(-\e i)^{b/(1+ab)}\ui^{(1+ab)^\e}(\e i)^{a/(1+ab)}$ 

for $i\in I$, $\e=\pm1$, $a,b$ in $K$;

(v) $\ui^a\ui^b=\ui^{ab}$, $\ui^{(1)}=1$ for $i\in I$, $a,b$ in $K$;

(vi) $\ui^a\uj^b=\uj^b\ui^a$ for $i,j$ in $I$, $a,b$ in $K$;

(vii) $\uj^a(\e i)^b=(\e i)^{a^{\e(i:j)}b}\uj^a$ for $i,j$ in $I$, $\e=\pm1$, $a,b$ in $K$;

(viii) $(\e i)^a(-\e j)^b=(-\e j)^b(\e i)^a$ for $i\ne j$ in $I$, $\e=\pm1$, $a,b$ in $K$.

\subhead 3.2\endsubhead
When $A=\left(\matrix 2&-1\\-1&2\endmatrix\right)$, $K=\RR_{>0}$, we can 
identify $G(K)$ with the submonoid of $SL_3(\RR)$ generated by:

$\left(\matrix 1&a&0\\0&1&0\\0&0&1\endmatrix\right)$, 
$\left(\matrix 1&0&0\\0&1&b\\0&0&1\endmatrix\right)$, 

$\left(\matrix 1&0&0\\c&1&0\\0&0&1\endmatrix\right)$, 
$\left(\matrix 1&0&0\\0&1&0\\0&d&1\endmatrix\right)$, 

$\left(\matrix e&0&0\\0&(1/e)&0\\0&0&1\endmatrix\right)$, 
$\left(\matrix 1&0&0\\0&f&0\\0&0&(1/f)\endmatrix\right)$, 

with $a,b,c,d,e,f$ in $\RR_{>0}$.

\subhead 3.3\endsubhead
The assignment $i^a\m i^a$ (with $i\in I,a\in K$) defines a monoid isomorphism of $U^+(K)$ onto a submonoid of 
$G(K)$; when $K=\{1\}$, we denote by $w\in G(\{1\})$ the image of $w\in U(\{1\})$ under this imbedding. 
The assignment $i^a\m(-i)^a$ (with $i\in I,a\in K$) defines a monoid isomorphism of $U^+(K)$ onto a submonoid of 
$G(K)$; when $K=\{1\}$, we denote by $-w\in G(\{1\})$ the image of $w\in U(\{1\})$ under this imbedding. 
The map $W\T W\to G(\{1\})$, $(w,w')\m w(-w')$ is a bijection of sets (not of monoids).

\subhead 3.4\endsubhead
Tits has said that $W$ ought to be regarded as the Chevalley 
group $G(\kk)$ where $\kk$ is the (non-existent) field with one element. 
But $G(\{1\})$ is defined for the semifield $\{1\}$. The bijections $W@>\si>>U^+(\{1\})$,
$W\T W@>\si>> G(\{1\})$ almost realizes the wish of Tits. 

\subhead 3.5\endsubhead
For general $K$, the semifield homomorphism $K\to\{1\}$ induces a
monoid homomorphism $G(K)\to G(\{1\})$. Let $G_{w,-w'}(K)$ be the 
fibre over $w(-w')$ of this homomorphism. We have $G(K)=\sqc_{(w,w')\in W\T W}G_{w,-w'}(K)$.
We now fix $(w,w')\in W\T W$. Let $M=|w|+|w'|+r$. 
It turns 
out that the set $G_{w,-w'}(K)$ can be parametrized by $K^M$, 
in fact in many ways, indexed by a certain finite set $\co_{w,-w'}$.
Let $\co_{-w'}$ be the set of sequences $(-i_1,-i_2,\do,-i_{|w'|})$ in $-I$ such 
that $(i_1,i_2,\do,i_{|w'|})\in\co_{w'}$. Let $\co_{w,-w'}$ be the set of sequences 
$(h_1,h_2,\do,h_M)$ in $I\sqc(-I)\sqc\uI$ such that the 
subsequence consisting of symbols in $I$ is in $\co_w$, the subsequence 
consisting of symbols in $-I$ is in $\co_{-w'}$, the subsequence consisting 
of symbols in $\uI$ contains each symbol $\ui$ (with $i\in I$) exactly once.

For $\hh=(h_1,h_2,\do,h_M)\in\co_{w,-w'}$ we define $\ps_\hh:K^M\to G_{w,-w'}(K)$ by 
$$\ps_\hh(a_1,a_2,\do,a_M)=h_1^{a_1}h_2^{a_2}\do h_M^{a_M}.$$
 This is a bijection.
The bijections $\ps_\hh:K^M\to G_{w,-w'}(K)$ (for various 
$\hh\in\co_{w,-w'}$) define a positive structure on $G_{w,-w'}(K)$.

In the case where $K=\RR_{>0}$ or $K=\RR(t)_{>0}$, the statements above 
are proved by using Bruhat decomposition in the group $G(\kk)$ considered in 2.1 with
$\kk=\RR$ or $\RR(t)$. (When $W$ is finite this is done in \cite{L19}. See also the
proof of \cite{L94, Lemma 2.3} and \cite{L94, 2.7}.) 
The case where $K=\ZZ$ follows from the case where $K=\RR(t)_{>0}$, using 
$\a:\RR(t)_{>0}@>>>\ZZ$ in 1.2. 

\head 4. Chevalley groups\endhead
\subhead 4.1\endsubhead
In this section we assume that $K=\RR_{>0}$ and that $I\ne\emp$. 
Let $\kk_0$ be a field and let $\kk$ be an algebraic closure of $\kk_0$. 

Let $w\in W,w'\in W$. Let $M=|w|+|w'|+r$. For $\hh,\hh'$ in $\co_{w,-w'}$, 
$\ps_{\hh'}\i\ps_\hh:K^M\to K^M$ (see 3.5) is of the form $(a_1,a_2,\do,a_M)\m(a'_1,a'_2,\do,a'_M)$ where 
$a'_s=(P^{\hh'}_\hh)_s(a_1,a_2,\do,a_M)/(Q_\hh^{\hh'})_s(a_1,a_2,\do,a_M)$ and each of 
$(P^{\hh'}_\hh)_s,(Q_\hh^{\hh'})_s$ is a nonzero polynomial in $\NN[X_1,X_2,\do,X_M]$ 
(independent of $K$) such that the g.c.d. of its $\ne0$ coeff. is $1$. 

Applying the obvious ring homomorphism $\ZZ\to\kk_0$ to the 
coefficients of these polynomials we obtain $\ne0$ polynomials $(\bP^{\hh'}_\hh)_s$,
$(\bQ_\hh^{\hh'})_s$ in $\kk_0[X_1,X_2,\do,X_M]$.
 We define a rational map $\bar\ps_\hh^{\hh'}:\kk^M\to\kk^M$ by 

$(z_1,z_2,\do,z_M)\m(z'_1,z'_2,\do,z'_m)$,

$z'_s=(\bP^{\hh'}_\hh)_s(z_1,z_2,\do,z_M)/(\bQ_\hh^{\hh'})_s(z_1,z_2,\do,z_M)$
\nl
This is a birational isomorphism. 
It induces an automorphism $[\ps_\hh^{\hh'}]$ of the quotient field 
$[\kk^M]$ of the coordinate ring of $\kk^M$. 
We have $[\ps_\hh^{\hh'}][\ps_{\hh'}^{\hh''}]=[\ps_\hh^{\hh''}]$ for any $\hh,\hh',\hh''$. 
Hence there is a well defined field $[G_{w,-w'}(\kk)]$ containing $\kk$ with 
$\kk$-linear field isomorphisms $[\ps_\hh]:[G_{w,-w'}(\kk)]\to[\kk^M]$ (for $\hh\in\co_{w,-w'}$) such 
that 

$[\ps_\hh^{\hh'}]=[\ps_\hh][\ps_{\hh'}]\i:[\kk^M]\to[\kk^M]$ for all $\hh,\hh'$.

\subhead 4.2\endsubhead
We now assume that $W$ is finite. Let $w_I,\nu$ be as in 2.6. Let $M=2\nu+r$. 
Let $i\in I,\e=\pm1,z\in\kk_0$. We can choose 
$\hh=(h_1,h_2,\do,h_M)\in\co_{\o,-\o}$ such that $h_1=\e i$. The isomorphism $\kk^M\to\kk^M$, 
$(z_1,z_2,\do,z_M)\m(z_1-z,z_2,\do,z_m)$ induces a field isomorphism $\t_z:[\kk^M]\to[\kk^M]$. 
Let $\bold A$ be the group of all $\kk$-linear field automorphisms of $[G_{\o,-\o}(\kk)]$. We 
define $(\e i)^z\in\bold A$ as the composition 
$$[G_{\o,-\o}(\kk)]@>[\ps_\hh]>>[\kk^M]@>\t_z>>[\kk^M]@>[\ps_\hh]\i>>[G_{\o,-\o}(\kk)].$$
Now $(\e i)^z$ is independent of the choice of $\hh$. Let $G(\kk_0)$ be the subgroup of $\bold A$ 
generated by $(\e i)^z$ for various $i\in I,\e=\pm1,z\in\kk_0$. Then $G(\kk_0)$ is the Chevalley group 
associated to $\kk_0$ and our Cartan matrix.

\head 5. Flag manifolds\endhead
\subhead 5.1\endsubhead
In this section $W$ is not assumed to be finite. We assume that $K$ is $\RR_{>0}$.
Let $G(\RR)$ be the group considered in 2.1.
Let $V$ be an $\RR$-vector space which is an irreducible highest weight integrable representation of $G(\RR)$
whose highest weight takes the value $1$ at any simple coroot.
Let $\et$ be a highest weight vector of $V$. Let $\BB$ 
be the canonical basis of $V$ (see \cite{L93, 11.10}) containing $\et$. Let $P$ be the set of lines in the
$\RR$-vector space $V$. Let $P_{\ge0}$ be the set of all $L\in P$ such that for some $x\in L-\{0\}$
all coordinates of $x$ with respect to the basis $\BB$ are $\ge0$. The flag manifold $\cb$ of $G(\RR)$  is defined 
as the subset of $P$ consisting of lines in the $G(\RR)$-orbit of the line spanned by $\et$. We define 
$\cb(K)=\cb\cap P_{\ge0}$. By a positivity property \cite{L93, 22.1.7} 
of $\BB$ (stated in the simply laced case but whose proof remains valid in our case), the
obvious $G(\RR)$-action on $\cb$ restricts to a $G(K)$-action on $\cb(K)$. (As mentioned in
2.1, $G(K)$ can be viewed as a submonoid of $G(\RR)$.)
When $W$ is finite, $\cb(K)$ is the same as the subset $\cb_{\ge0}$ defined in \cite{L94,\S8}. (This
follows from \cite{L94, 8.17}.)

\enddocument